\documentclass[12pt]{amsart}

\usepackage[all]{xy}
\usepackage{a4wide}
\usepackage{geometry}
\usepackage{amsmath}
\usepackage{amssymb}
\usepackage{amsthm}
\usepackage{amscd}
\usepackage{subfigure}
\usepackage{epsfig}
\usepackage{color}
\usepackage{graphicx}   
\usepackage{stmaryrd}   
\usepackage{mathrsfs}

\usepackage{changes}
\usepackage[skip=2pt,font=scriptsize]{caption}

\newcommand{\T}{\mathbb T}
\newcommand{\N}{\mathbb N}
\newcommand{\bee}{\begin{equation}}
\newcommand{\eee}{\end{equation}}

\newcommand{\Lb}{\mbox {\boldmath ${\Lambda}$}}
\newcommand{\Gb}{\mbox {\boldmath ${\Gamma}$}}

\def\ov{\overline}

\def\b0{{\bf 0}}

\definecolor{darkgreen}{rgb}{0.2,0.5,0}

\newcommand{\be}{\begin{eqnarray}}
\newcommand{\ee}{\end{eqnarray}}
\newcommand{\supp}{\mbox{\rm supp}}

\newcommand{\R}{{\mathbb R}}
\newcommand{\Q}{{\mathbb Q}}
\newcommand{\Z}{{\mathbb Z}}
\newcommand{\C}{{\mathbb C}}

\newcommand{\Ak}{{\mathcal A}}

\newcommand{\Sf}{{\mbox{\sf S}}}
\newcommand{\Dk}{{\mathcal D}}

\newcommand{\Pk}{{\mathcal P}}

\newcommand{\Tk}{{\mathcal T}}
\newcommand{\Xk}{{\mathcal X}}

\newcommand{\bx}{{\bf x}}

\newcommand{\Lam}{{\Lambda}}

\newcommand{\lam}{\lambda}

\newcommand{\om}{\omega}

\newcommand{\balpha}{\mbox{\boldmath{$\alpha$}}}
\newcommand{\bbeta}{\mbox{\boldmath{$\beta$}}}

\newcommand{\bgamma}{\mbox{\boldmath{$\gamma$}}}

\newtheorem{theorem}{Theorem}[section]
\newtheorem{lemma}[theorem]{Lemma}

\newtheorem{prop}[theorem]{Proposition}
\newtheorem{cor}[theorem]{Corollary}

\newtheorem{exam}[theorem]{Example}
\newtheorem{remark}[theorem]{Remark}

\newtheorem{defi}[theorem]{\it{Definition}}
\numberwithin{equation}{section}

\begin{document}

\title[Pure discrete spectrum and regular model sets]{Pure discrete spectrum and regular model sets in unimodular substitution tilings on $\R^d$}







\maketitle

\centerline{Dong-il Lee$^{1}$, Shigeki Akiyama $^{2}$ and   Jeong-Yup Lee  $^{3}$*}

\hspace*{4em} 

\noindent
{\footnotesize
$^{1}$ \quad   Department of Mathematics, Seoul Women's University, Seoul 01797, Korea \\
$^{2}$ \quad 
Institute of Mathematics, University of Tsukuba, 1-1-1 Tennodai, Tsukuba, Ibaraki, Japan \\
\mbox{~} \quad  (zip:305-8571); akiyama@math.tsukuba.ac.jp \\
$^{3}$ \quad
 Department of Mathematics Education, Catholic Kwandong University,
 Gangneung, Gangwon \\
\mbox{~}  \quad 210-701, Korea ; jylee@cku.ac.kr, \ \  KIAS, 85 Hoegiro, Dongdaemun-gu, Seoul, 02455, Korea
}



\maketitle

\begin{abstract}
We consider primitive substitution tilings on $\R^d$ whose expansion maps are unimodular. We assume that all the eigenvalues of the expansion maps are algebraic conjugates with the same multiplicity. In this case, we can construct a cut-and-project scheme with a Euclidean internal space. 
Under some additional condition, we show that if the substitution tiling has pure discrete spectrum, then the corresponding representative point sets are regular model sets in that cut-and-project scheme.  

\hspace*{4em}
\vspace{5mm}


\end{abstract}

\section{Introduction}

In the study of aperiodic tilings, it has been an interesting problem to characterize pure discrete spectrum of tiling dynamical systems. 
 This property is related to understanding the structure of mathematical quasicrystals. 
 For this direction of study, substitution tilings have been good models, since they have highly symmetrical structures.
A lot of research has been done in this direction (see (Akiyama \& Barge \& Berthe \& Lee \& Siegel 2015), (Baake \& Grimm 2013) and references therein).  
Given a substitution tiling with pure discrete spectrum, it is known that this can be described via cut-and-project scheme(CPS) (Lee 2007).
However, in (Lee 2007), the construction of the cut-and-project scheme is with an abstract internal space built from the pure discrete spectral property. Since the internal space is an abstract space, it is neither easy to understand the tiling structure, nor is clear if the model sets are regular or not.
In the case of $1$-dimensional substitution tilings with pure discrete spectrum, it is shown that CPS with an Euclidean internal space
can be built and the corresponding representative point sets are regular model sets in (Barge \& Kwapisz 2006). In this paper, we consider 
substitution tilings on $\R^d$ with pure discrete spectrum whose expansion maps are unimodular. We show that it is possible to construct a CPS with an Euclidean internal space and show that the corresponding representative point sets are regular model sets in that CPS.

The outline of the paper is as follows. First, we consider a repetitive primitive substitution tiling on $\R^d$ whose expansion map is unimodular.
Then we build a CPS with an Euclidean internal space in section \ref{CPS}. In section \ref{KnownResults}, we discuss some known results around pure discrete spectrum, Meyer set, and Pisot family. In section \ref{MainResult}, under the assumption of pure discrete spectrum, we show that each representative point set of the tiling is actually a regular model set in the CPS with the Euclidean internal space.

\section{Preliminaries}

\subsection{Tilings}

We begin with a set of types (or colours) $\{1,\ldots, \kappa \}$, which
we fix once and for all. A {\em tile} in $\R^d$ is defined as a
pair $T=(A,i)$ where $A=\supp(T)$ (the support of $T$) is a
compact set in $\R^d$, which is the closure of its interior, and
$i=l(T)\in \{1,\ldots, \kappa \}$ is the type of $T$. 
We let $g+T =
(g+A,i)$ for $g\in \R^d$. 
We say that a set $P$ of tiles is a {\em
patch} if the number of tiles in $P$ is finite and the tiles of
$P$ have mutually disjoint interiors. The {\em support of a patch}
is the union of the supports of the tiles that are in it. The {\em
translate of a patch} $P$ by $g\in \R^d$ is $g+P := \{g+T:\ T\in
P\}$. We say that two patches $P_1$ and $P_2$ are {\em
translationally equivalent} if $P_2 = g+P_1$ for some $g\in \R^d$.
A {\em tiling} of $\R^d$ is a set $\Tk$ of tiles such that $\R^d =
\bigcup \{\supp(T) : T \in \Tk\}$ and distinct tiles have disjoint
interiors. 
We
always assume that any two $\Tk$-tiles with the same colour are
translationally equivalent (hence there are finitely many
$\Tk$-tiles up to translations). 
Given a tiling $\Tk$, a finite set of tiles of $\Tk$ is called {\em $\Tk$-patch}.
Recall that a tiling $\Tk$ is said to be {\em repetitive} if every $\Tk$-patch occurs relatively dense in space, up to translation.
We say that a tiling $\Tk$ has {\em finite local complexity(FLC)} if for every $R > 0$, there are finitely many equivalence classes of $\Tk$-patches of diameter less than $R$. 

\subsection{Delone $\kappa$-sets}

A {\em $\kappa$-set} in $\R^d$ is a
subset $\Lb = \Lam_1 \times \dots \times \Lam_{\kappa}
\subset \R^d \times \dots \times \R^d$ \; ($\kappa$ copies)
where $\Lam_i \subset \R^d$ and $\kappa$ is the number of colours. 
We also write
$\Lb = (\Lam_1, \dots, \Lam_{\kappa}) = (\Lam_i)_{i\le \kappa}$.
Recall that a Delone set is a relatively dense and uniformly discrete subset of $\R^d$.
We say that $\Lb=(\Lambda_i)_{i\le \kappa}$ is a {\em Delone $\kappa$-set} in $\R^d$ if
each $\Lambda_i$ is Delone and $\supp(\Lb):=\bigcup_{i=1}^{\kappa} \Lambda_i \subset \R^d$ is Delone. 
The types (or colours) of points for Delone $\kappa$-sets have the meaning analogous
to the colours of tiles for tilings. 
We define repetitivity and FLC for Delone $\kappa$-set in the same way as for tilings.
A Delone set $\Lambda$ is called a {\em Meyer set} in $\R^d$ if $\Lambda - \Lambda$ is uniformly discrete, which is equivalent to say that $\Lambda -\Lambda = \Lambda + F$ for some finite set $F$ (see (Moody 1997)). If $\Lb$ is a Delone $\kappa$-set and $\supp(\Lb)$ is a Meyer set, we say that $\Lb$ is a Meyer set.

\subsection{Substitutions}

We say that a linear map $\phi:\ \R^d \to \R^d$ is {\em expansive} if there is a constant
$c > 1$ with 
\[ d(\phi x, \phi y) \ge c \cdot d(x, y) \]
for all $x, y \in \R^d$ under some metric $d$ on $\R^d$ compatible with the standard topology.

\begin{defi}\label{def-subst}
{\em Let $\Ak = \{T_1,\ldots,T_{\kappa} \}$ be a finite set of tiles on $\R^d$
such that $T_i=(A_i,i)$; we will call them {\em prototiles}.
Denote by $\Pk_{\Ak}$ the set of
patches made of tiles each of which is a translate of one of $T_i$'s.
We say that $\omega: \Ak \to \Pk_{\Ak}$ is a {\em tile-substitution} (or simply
{\em substitution}) with an expansive map $\phi$ if there exist finite sets $\Dk_{ij}\subset \R^d$ for
$i,j \le \kappa$, such that
\begin{equation}
\om(T_j)=
\{u+T_i:\ u\in \Dk_{ij},\ i=1,\ldots,\kappa \}
\label{subdiv}
\end{equation}
with
\be \label{eq-til} 
\phi A_j = \bigcup_{i=1}^{\kappa} (\Dk_{ij}+A_i) \ \ \  \mbox{for each} \  j \le \kappa.
\ee
Here all sets in the right-hand side must have disjoint interiors;
it is possible for some of the $\Dk_{ij}$ to be empty.  We call the finite set $\Dk_{ij}$ a {\em digit set} (Lagarias \& Wang 1996). The {\em substitution $\kappa\times\kappa$ matrix} $\Sf$ of the tile-substitution is defined by $\Sf(i,j) = \# \Dk_{ij}$. We say that $\phi$ is {\em unimodular} if 
the minimal polynomial of $\phi$ over $\Q$ has constant term $\pm 1$(i.e. $\mbox{det} \ \phi = \pm 1 $), that is to say, the product of all roots of the minimal polynomial of $\phi$ is $\pm 1$. 
}
\end{defi}

Note that for $M \in \N$,
\[ \phi^M A_j = \bigcup_{i=1}^{\kappa} (\Dk^M_{ij}+A_i) \ \ \  \mbox{for} \  j \le \kappa,   \]
where
\begin{eqnarray}  \label{M-th-iteration-digitSet}
  (\mathcal{D}^M)_{ij} = \bigcup_{k_1, k_2, \dots, k_{(M-1)} \le \kappa} (\mathcal{D}_{ik_1} + \phi \mathcal{D}_{k_1 k_2} +
\cdots + \phi^{M-1} \mathcal{D}_{k_{(M-1)}j})
\end{eqnarray}

The tile-substitution is extended to translated prototiles by 
\be \label{transl}
\om(T_j - x) = \om(T_j) - \phi x.
\ee 
The equations (\ref{eq-til}) allow us to extend $\om$ to patches in $\Pk_\Ak$ defining by
$\om(P) = \bigcup_{T\in P} \om(T)$. It is similarly extended to tilings all of whose tiles are translates of the prototiles from $\Ak$. A tiling $\Tk$ satisfying $\om(\Tk) = \Tk$ is called a
{\em fixed point of the tile-substitution}, or a {\em substitution tiling with expansion map $\phi$}. 
It is known that one can always find a periodic point for $\om$ in the tiling dynamical hull, i.e. $\om^N(\Tk) = \Tk$ for some $N \in \N$.  In this case we use $\om^N$ in the place of $\om$ to obtain a  fixed point tiling.     
We say that the substitution tiling $\Tk$ is {\em primitive}, if there is an $\ell > 0$ for which $\Sf^{\ell}$ has no zero entries, where $\Sf$ is the substitution matrix. 

\begin{defi} \label{def-subst-mul}
{\em $\Lb = (\Lam_i)_{i\le \kappa}$ is called a {\em
substitution Delone $\kappa$-set} if $\Lb$ is a Delone $\kappa$-set and
there exist an expansive map
$\phi:\, \R^d\to \R^d$ and finite sets $\Dk_{ij}$ for $i,j\le \kappa $ such that
\be \label{eq-sub}
\Lambda_i = \bigcup_{j=1}^{\kappa} (\phi \Lambda_j + \Dk_{ij}),\ \ \ i \le \kappa,
\ee
where the unions on the right-hand side are disjoint. 
}
\end{defi}

There is a standard way to choose distinguished points in the tiles of primitive substitution
tiling so that they form a $\phi$-invariant Delone $\kappa$-set. They are called {\em control points} (Thurston 1989, Praggastis 1999) which are defined below.

\begin{defi}  
{\em Let $\Tk$ be a fixed point of a primitive substitution with an
expansion map $\phi$. For every $\Tk$-tile $T$, we choose a tile
$\gamma T$ on the patch $\omega(T)$;  for all tiles of the same type, we choose
$\gamma T$ with the same relative position. This defines a map
$\gamma: \Tk \to \Tk $ called the {\em tile map}. Then we define the {\em control point} for a tile $T \in \Tk $ by
\begin{eqnarray}  \label{tileMapDef}
\{ c(T)\} = \bigcap_{m=0}^{\infty} \phi^{-m} (\gamma^m T).
\end{eqnarray}
}
\end{defi}

\medskip

The control points satisfy the following:
\begin{itemize}
\item[(a)] $T' = T + c(T') - c(T)$, for any tiles $T, T'$ of the same type;
\item[(b)] $\phi(c(T)) = c(\gamma T)$, for $T \in \Tk $.
\end{itemize}
For tiles of any tiling $\mathcal{S} \in  X_{\Tk}$, the control points have the same
relative position as in $\Tk$-tiles. The choice of control points is non-unique,
but there are only finitely many possibilities, determined by the choice of the tile
map. 
Let
\begin{eqnarray}  \label{translation-vectorSet}
  \Xi = \bigcup_{i=1}^{\kappa} (\mathcal{C}_i - \mathcal{C}_i). 
\end{eqnarray}  
It is possible to consider a tile map 
\begin{eqnarray} \label{the-tile-map}
\gamma: \Tk \to \Tk \  \mbox{ s.t  $\forall$ \ $T \in \Tk$, the tile $\gamma(T)$ has the same tile type in $\Tk$}. 
\end{eqnarray}
Then for any $T, S \in \Tk$, 
\[ c(\gamma T) - c(\gamma S) \ \in \  \Xi. \] 
Let
\[ \mathcal{C}:= \mathcal{C}(\Tk) = (\mathcal{C}_i)_{i \le \kappa} = \{ c(T): T \in \Tk\}  \]
be a set of control points of the tiling $\Tk$ in $\R^d$. In what follows, if there is no confusion, we will use the same notation $\mathcal{C}$ to mean $\supp \ \mathcal{C}$. 
For the main results of this paper, we need the property that 
$\langle \bigcup_{i\le \kappa} \mathcal{C}_i \rangle = \langle \Xi \rangle$
with $0\in \bigcup_{i\le \kappa} \mathcal{C}_i$.
Under the assumption that $\phi$ is unimodular, 
this can be
achieved by taking a proper control point set which comes from a certain tile map.
We define the tile map as follows.
It is known that there exists a finite patch $\mathcal{P}$ in a primitive substitution tiling which generates the whole tiling $\Tk$ (Lagarias-Wang 2003). Although it was defined with primitive substitution point sets in (Lagarias-Wang 2003), it is easy to see that the same property holds for primitive substitution tilings. We call the finite patch $\mathcal{P}$ {\em generating tile set}.
When we apply the substitution infinitely many times to the generating tile set $\mathcal{P}$, we obtain the whole substitution tiling. So there exists $n \in \N$ such that $n$-th iteration of the substitution to the generating tile set covers the origin. 
We choose a tile $R$ in a patch $\omega^n (\mathcal{P})$ which contains the origin, where $R = \alpha + T_j$ for some $1 \le j \le \kappa$. 
Then there exists a fixed tile $S \in \mathcal{P}$ such that $R \in \omega^n(S)$. Replacing the substitution $\omega$ by $\omega^n$, we can define a tile map $\gamma$ so that 
\[  \left\{  
\begin{array}{ll}
\mbox{ $\gamma(T)$ is a $j$-type tile in $\omega^n(T)$} \ \ \ \mbox{if $T \in \mathcal{T}$ with $T \neq S $} \\
\gamma(S) = R.
 \end{array}
   \right.
\]
Then $0 \in \mathcal{C}_j$ by the definition of the control point sets and so $\mathcal{C}_j \subset \Xi$.
Notice that 
\begin{eqnarray} \label{phi-inverse-unimodularity}
\phi^{-1}\in \Z[\phi],
\end{eqnarray} 
since $\phi$ is unimodular. 
From the construction of the tile map, we have 
$\phi(\mathcal{C}_i) \subset \mathcal{C}_j$ for any $1 \le i \le \kappa$. From (\ref{phi-inverse-unimodularity}), 
 $\mathcal{C}_i \subset \langle \mathcal{C}_j \rangle$ for any $1 \le i \le \kappa$. Hence
$\bigcup_{i \le \kappa} \mathcal{C}_i \subset \langle \mathcal{C}_j \rangle$.
Thus
\begin{eqnarray} \label{control-point-set-translation-set}
\left\langle \bigcup_{i\le \kappa} \mathcal{C}_i \right\rangle = \langle \Xi \rangle.
\end{eqnarray}

\begin{remark} 
{\em In the case of primitive unimodular irreducible $1$-dimension Pisot substitution tilings, it is known that 
$\langle \bigcup_{i\le \kappa} \mathcal{C}_i \rangle = \langle \Xi \rangle. $
 See (Barge \& Kwapisz 2006, Sing 2006).
}

\end{remark}

\subsection{Pure point spectrum and algebraic coincidence}

Let $X_{\Tk}$ be the collection of tilings on $\R^d$ each of whose patches is a translate of a $\Tk$-patch. In the case that $\Tk$ has FLC, we give a usual metric $\delta$ on the tilings in such a way that two tilings are close if there is a large agreement on a big region with small shift (see (Schlottmann 2000, Radin \& Wolff 1992, Lee \& Moody \& Solomyak 2003)). Then 
\[ X_{\Tk}= \overline{\{-h+\Tk : h \in \R^d \}} \] 
where the closure is taken in the topology induced by the metric $\delta$. 
For non-FLC tilings, we can consider ``local rubber topology'' on the collection of tilings (M\"{u}ller \& Richard 2013, Lenz \& Stollmann 2003, Baake \& Lenz 2004, Lee \& Solomyak 2019) and obtain $X_{\Tk}$ as the completion of the orbit closure of $\Tk$ under this topology. For tilings with FLC, the two topologies coincide. In the either case of FLC or non-FLC tilings, we obtain a compact space $X_{\Tk}$.
We have a natural action of $\R^d$ on the dynamical hull $X_{\Tk}$ of $\Tk$ by translations and get a topological dynamical system $(X_{\Tk}, \R^d)$.
Let us assume that there is a unique ergodic measure $\mu$ in the dynamical system $(X_{\Tk}, \R^d)$ and consider the measure preserving dynamical system $(X_{\Tk}, \mu, \R^d)$. It is known that a dynamical system $(X_{\Tk}, \R^d)$ with a primitive substitution tiling $\Tk$ always has a unique ergodic measure in (Solomyak 1997,  Lee \& Moody \& Solomyak 2003).
We consider the associated group of unitary operators $\{T_x\}_{x\in \R^d}$ on
$L^2(X_{\Tk},\mu)$:
$$
T_x g(\Tk') = g(-x + \Tk').
$$
Every $g \in L^2(X_{\Tk},\mu)$ defines a function on $\R^d$ by
$x \mapsto \langle T_x g,g \rangle$.
This function is positive definite on $\R^d$, so its
Fourier transform is a positive measure $\sigma_g$ on $\R^d$ called the
{\em spectral measure} corresponding to $g$.
The dynamical system $(X_{\Tk}, \mu, \R^d)$ is said to have {\em
pure discrete spectrum} if $\sigma_g$ is pure point for every $g \in
L^2(X_{\Tk}, \mu)$. We also say that $\Tk$ has pure
discrete spectrum if the dynamical system $(X_{\Tk}, \mu, \R^d)$ has
pure discrete spectrum.

When we restrict to primitive substitution tilings, we note that a tiling $\Tk$ has pure discrete spectrum if and only if the control point set $\mathcal{C}(\Tk)$ of the tiling $\Tk$ admits an algebraic coincidence (See Definition \ref{algebraic-coincidence}). So from now on when we assume pure discrete spectrum for $\Tk$, we can directly use the property of algebraic coincidence. We give the corresponding definition and theorem below.

\begin{defi} \label{algebraic-coincidence}
{\em 
Let $\Tk$ be a primitive substitution tiling on $\R^d$ with an expansive map $\phi$ and $\mathcal{C} = (\mathcal{C}_i)_{i \le \kappa}$ be a corresponding control point set. We say that $\mathcal{C}$ admits an {\em algebraic coincidence} if there exists $M \in \Z_+$ and $\xi \in \mathcal{C}_i$ for some $1 \le i \le \kappa$ such that 
\[ \xi + \phi^M \Xi \ \subset \ \mathcal{C}_i. \] 
Here recall from (\ref{translation-vectorSet}) that
$ \Xi = \bigcup_{i=1}^{\kappa} (\mathcal{C}_i - \mathcal{C}_i)$.
}
\end{defi}

\medskip

Note that if the algebraic coincidence is assumed, then for some $M \in \Z_+$,
\begin{eqnarray}  \label{translation-vectorSet-relation} 
\phi^M \Xi - \phi^M \Xi \ \subset \ \mathcal{C}_i - \mathcal{C}_i \  \subset \ \Xi.
\end{eqnarray}

\begin{theorem} [Theorem 3.13 (Lee 2007)] \label{pdd-algCoin}
Let $\Tk$ be a primitive substitution tiling on $\R^d$ with an expansive map $\phi$ and $\mathcal{C} = (\mathcal{C}_i)_{i \le \kappa}$ be a corresponding control point set. Suppose that all the eigenvalues of $\phi$ are algebraic integers. 
Then $\Tk$ has pure discrete spectrum if and only if $\mathcal{C}$ admits an algebraic coincidence.
\end{theorem}

The above theorem is stated with FLC in (Lee 2007). But from Lemma \ref{pdd-imply-relDenseEV} and Proposition \ref{relDense-Meyer}, pure discrete dynamical spectrum of $\Tk$ implies the Meyer property of the control point set $\mathcal{C}$.  All Meyer sets have FLC. So it is a consequence of pure discrete dynamical spectrum.
On the other hand, the algebraic coincidence implies that  
\[ \xi + \phi^M \Xi \ \subset \ \mathcal{C}_i \ \ \ \mbox{for some $\xi \in \mathcal{C}_i$ and $1 \le i \le \kappa$}. \]
This means that $\phi^M \Xi$ is uniformly discrete and thus $\Xi$ is uniformly discrete. From (Moody 1997), we obtain that $\Xi - \Xi$ is uniformly discrete.
For any $1 \le i, j \le \kappa$,
\[ \mathcal{C}_i - \mathcal{C}_j \subset a+ \Xi - (b+ \Xi) = a - b 
+ \Xi - \Xi \ \ \ \mbox{for some $a \in \mathcal{C}_i$ and $b \in \mathcal{C}_j$}. \]
Hence $\mbox{supp}(\mathcal{C}) = \bigcup_{i=1}^{\kappa} \mathcal{C}_i$ is a Meyer set (Moody 1997). 
So it is not necessary to assume FLC here.

\subsection{Cut-and-project scheme}

We give definitions for cut and project scheme and model sets constructed with $\R^d$ and a locally compact Abelian group.  
\begin{defi} \label{cut-and-project1}
{\em A {\em {cut and project scheme}} (CPS) consists of a collection of spaces and mappings as follows;
\be
\begin{array}{ccccc} \label{cut-and-project-scheme1}
 \R^d & \stackrel{\pi_{1}}{\longleftarrow} & \R^d \times H & \stackrel{\pi_{2}}
{\longrightarrow} & H \\
 && \bigcup \\
 && \widetilde{L}
\end{array}
\ee
where $\R^d$ is a real Euclidean space, $H$ is a locally
compact Abelian group, $\pi_{1}$ and $\pi_{2}$ are the canonical projections,
$ \widetilde{L} \subset {\R^d \times H}$ is a lattice,  i.e.\
a discrete subgroup for which the quotient group
$(\R^d \times H) / \widetilde{L}$ is
compact, $\pi_{1}|_{ \widetilde{L}}$ is injective,
and $\pi_{2}(\widetilde{L})$ is dense in $H$.
For a subset $V \subset H$, we denote 
\[\Lam(V) := \{ \pi_{1}(x) \in \R^d : x \in \widetilde{L}, \pi_{2}(x) \in V\}.\]
A {\em model set} in $\R^d$ is a subset $\Lam$ of $\R^d$ of the form $\Lam(W)$, 
where $W \subset H$ has non-empty interior and compact closure.
The model set $\Lam$ is {\em regular} if the boundary of $W$
\[ \partial{W} = \overline{W} \ \backslash \ W^{\circ}
\] 
is of (Haar) measure $0$.
We say that $\Lb = (\Lam_i)_{i \le \kappa}$ is a {\em model $\kappa$-set} (resp. {\em regular model $\kappa$-set}) if each $\Lam_i$ is a model set (resp. regular model set) with respect to the same CPS. Especially when $H$ is an Euclidean space, we call the model set $\Lambda$ a {\em Euclidean model set} (see (Baake \& Grimm 2013)).}
\end{defi}

\section{Cut-and-project scheme on substitution tilings} \label{CPS}

Throughout the rest of the paper, we assume that $\phi$ is diagonalizable, the eigenvalues of $\phi$ are algebraically conjugate with the same multiplicity, since 
the structure of a module generated by the control points is known only under this assumption (Lee \& Solomyak 2012).

Let 
\begin{eqnarray} \label{real-eigenvalues}
\lam'_{1}, \ldots, \lam'_{s}
\end{eqnarray} be the distinct real eigenvalues of $\phi$ and 
\begin{eqnarray} \label{complex-eigenvalues}
\lam'_{s+1},\ov{\lam'}_{s+1},\ldots,\lam'_{s + t},\ov{\lam'}_{s+t}
\end{eqnarray}
  be the distinct complex eigenvalues of $\phi$. By the above assumption, all these eigenvalues appear with the same multiplicity, which we will denote by $J$.
  Recall that $\phi$ is assumed to be diagonalizable over $\C$. For a complex eigenvalue $\lambda$ of $\phi$, the $2 \times 2$ diagonal block 
  $\left[ \begin{array}{cc}
  \lambda & 0 \\
  0 & \overline{\lambda}
  \end{array}
  \right]
  $ is similar to a real $2 \times 2$ matrix 
  \begin{eqnarray} \label{transition-of-matrix}
   \left[ 
  \begin{array}{cc}
  a & -b \\
  b & a
  \end{array} \right] =  
  S^{-1} \left[ \begin{array}{cc}
  \lambda & 0 \\
  0 & \overline{\lambda}
  \end{array} \right] S, 
  \end{eqnarray}
  where $\lambda = a+ ib, a, b \in \R $, and 
  $S= \frac{1}{\sqrt{2}} \left[ \begin{array}{cc}
  1 & i \\
  1 & -i
  \end{array} \right].$
Since $\phi$ is diagonalizable, by eventually changing the basis in $\R^d$, we can assume without loss of generality that
\begin{eqnarray}  \label{matrix-form}
\phi = {\small \left[
\begin{array}{ccc}
                  \psi_1 & \cdots & {O} \\
                  \vdots & \ddots & \vdots \\
                  {O} & \cdots & \psi_J
                  \end{array} \right] } \  \mbox{and} \ 
                  \psi_j = \psi :=  {\small \left[
                  \begin{array}{ccc}
                  A_{1} & \cdots & 0 \\
                  \vdots & \ddots & \vdots \\
                  0 & \cdots & A_{s+t}
                  \end{array} \right] } 
\end{eqnarray}
where $A_{k}$ is a real $1\times 1$ matrix for $1 \le k \le s$, a
real $2\times 2$ matrix of the form $\left[
\begin{array}{lr}
         a_{k} & -b_{k} \\
         b_{k} & a_{k}
        \end{array} \right] $
for $s + 1 \le k \le s + t$, and ${O}$ is the $(s+2t) \times (s+2t)$
zero matrix, and $1 \le j \le J$. 

Let 
\[ m = s +2t . \]
Note that $m$ is the degree of 
the minimal polynomial of $\phi$ over $\Q$.
For each $1 \le j \le J$, let
\begin{equation}
H_{j} = \{0\}^{(j-1)m} \times \R^{m} \times
\{0\}^{d- jm}\,. \nonumber
\end{equation}
Further, for each $H_{j}$ we have the direct sum decomposition
$$
H_j = \bigoplus_{k=1}^{s + t} E_{jk},
$$
such that each $E_{jk}$ is $\phi,\phi^{-1}$-invariant and
$\phi|_{E_{jk}} \approx A_{k}$, identifying $E_{jk}$ with $\R$ or
$\R^2$. 
Let 
\[ \phi_{j} = \phi|_{H_{j}}. \]

Let $P_{j}$ be the canonical projection of $\R^d$ onto $H_{j}$
such that \be \label{def-projection} P_{j}(\bx) = \bx_{j}, \ee
where $\bx = \bx_{1} + \cdots + \bx_{J}$ and $\bx_{j} \in H_{j}$
with $1 \le j \le J$.

We define $\balpha_{j} \in H_{j}$ such that for each $1 \le k \le
d$, \be \label{def-alpha}
 (\balpha_{j})_{k} =
\left\{\begin{array}{ll}
                             1 \ \ \ & \mbox{if} \ \  (j-1)m + 1 \le k \le jm;\\
                             0 \ \ \ & \mbox{else} \,.
                           \end{array} \right.
\ee

We recall the following theorem for the module structure of the control point sets.
From (Lemma 6.1 (Lee \& Solomyak 2012)), we can readily obtain 
the property:\footnote{
This fact is stated in a slightly different way in (Lee \& Solomyak 2012).}
\begin{eqnarray} \label{alphas-are-contorlPoints}
\balpha_1, \dots, \balpha_J \in \sigma(\ \mathcal{C}(\Tk))
\end{eqnarray}
which is used in
the proof of Lemma \ref{latticePoints-in-translationVectorSet}.
So we state (Theorem 4.1 (Lee \& Solomyak 2012)) in the following form. 
Let
$$
\Q[\phi]:= \{p(\phi):\ p(x) \in \Q[x]\},\ \ \Z[\phi]:= \{p(\phi):\ p(x) \in \Z[x]\}.
$$
\begin{theorem} [Theorem. 4.1 (Lee \& Solomyak 2012)] \label{controlPoints-Lattice}
Let $\Tk$ be a repetitive primitive substitution tiling on $\R^d$ with an expansion map $\phi$.
Assume that $\Tk$ has FLC, $\phi$ is diagonalizable, and all the eigenvalues of $\phi$ are algebraically conjugate with the same multiplicity $J$.
Then there exists an isomorphism $\sigma: \R^d \to \R^d $ such that
\[ \sigma \phi = \phi \sigma\ \ \ \mbox{and} \ \
 \sigma({ \mathcal{C}(\Tk)}) \subset 
\Z[\phi]{{\balpha}_1} + \cdots +  \Z[\phi]{{\balpha}_J} \,, 
\]
where ${\balpha}_j$, $1 \le j \le J$, are given as (\ref{def-alpha}), and $\balpha_1, \dots, \balpha_J \in \sigma({ \mathcal{C}(\Tk)}) $. 
\end{theorem}

\medskip

Since $\phi$ is a block diagonal matrix as shown in (\ref{matrix-form}),
we can note that $\balpha_1, \dots, \balpha_J$ are linearly independent over $\Z[\phi]$.
A tiling $\Tk$ is said to be {\em rigid} if $\Tk$ satisfies the result of Theorem \ref{controlPoints-Lattice}, that is to say,
there exists a linear isomorphism $\sigma: \R^d \to \R^d$ such that
\[ \sigma \phi = \phi \sigma\ \ \ \mbox{and} \ \
 \sigma({ \mathcal{C}(\Tk)}) \subset 
\Z[\phi]{{\balpha}_1} + \cdots +  \Z[\phi]{{\balpha}_J} \,, \]
where ${\balpha}_j$, $1 \le j \le J$, are given as (\ref{def-alpha}).

\medskip

\subsection{Construction of cut-and-project scheme}

Consider that $\phi$ is unimodular and diagonalizable, all the eigenvalues of $\phi$ are algebraic integers and
algebraically conjugate with the same multiplicity $J$, and $\Tk$ is rigid. 
Since $\phi$ is an expansion map and unimodular,  there exists at least one other algebraic conjugate other than eigenvalues of $\phi$.

It is known that if $\phi$ is a diagonalizable expansion map of a primitive substitution tiling with FLC, every eigenvalues of $\phi$ is an algebraic integer (Kenyon \& Solomyak 2010). So it is natural to assume that all the eigenvalues of $\phi$ are algebraic integers in the assumption.
In (\ref{matrix-form}), suppose that the minimal polynomial of $\psi$ over $\Q$ has
 $e$ number of real roots, and $f$ number of pairs of complex conjugate roots. Recall that 
 \[ \lam'_{1}, \ldots, \lam_{s} = \lam'_{s}, \ldots, \lam'_{s+1}, \ldots, \ov{\lam'_{s+t}} \] are distinct eigenvalues of $\phi$ from (\ref{real-eigenvalues}) and (\ref{}). 
 Let us consider the roots in the following order: 
\[
\lambda_1, \dots,  \lambda_s, \lambda_{s+1}, \dots. \lambda_{e}, \lambda_{e+1}, \overline{\lambda_{e+1}}, \dots,  \lambda_{e+t}, \overline{\lambda_{e+t}},  \lambda_{e+t+1}, \overline{\lambda_{e+t+1}}, \dots, \lambda_{e+f}, \overline{\lambda_{e+f}}, 
\]
for which
\[ \lam_{1}=\lam'_{1}, \ \ldots, \ \lam_{s} = \lam'_{s}, \] 
\[ \lam_{e+1}=\lam'_{s+1},\ov{\lam_{e+1}}= \ov{\lam'_{s+1}}, \ \ldots, \ \lam_{e+t}=\lam'_{s + t},  \ov{\lam_{e+t}}=\ov{\lam'_{s+t}}, \]
where $\lam'_{1}, \ldots, \lam'_{s}, \lam'_{s+1}, \ov{\lam'_{s+1}}, \ldots, \lam'_{e+t}, \ov{\lam'_{s+t}}$ are same as in (\ref{real-eigenvalues}) and  (\ref{complex-eigenvalues}).

 Let 
\begin{eqnarray} \label{n=e+2f}
n = e + 2f. 
\end{eqnarray}
We consider a space where the rest of roots of the minimal polynomial of $\psi$ other than the eigenvalues of $\psi$ lies. 
Using similar matrices as in (\ref{transition-of-matrix}) we can consider the space as an Euclidean space. 
Let 
\[ \mathbb{G}_j:= \R^{n-m}, \ \ \ 1 \le j \le J , \]

For $1 \le j \le J$, define a $(n-m) \times (n-m)$ matrix
\[ D_j := \left[ 
\begin{array}{cccccc}
 A_{s+1} & \cdots & 0 & 0 & \cdots & 0 \\
 \vdots & \ddots & \vdots & \vdots & \ddots & \vdots \\ 
  0 & \cdots &   A_{e} & 0 & \cdots & 0 \\ 
   0 & \cdots & 0 &  A_{e+t+1} & \cdots & 0 \\  
  \vdots & \ddots & \vdots & \vdots & \ddots & \vdots \\
   0 & \cdots & 0 &  0 & \cdots & A_{e+f} 
      \end{array}
\right]
\]
where $A_{s+ \imath}$ is a real $1 \times 1$ matrix with the value $\lambda_{s+ \imath}$ for $1 \le \imath \le e-s$, and $A_{e+t+\jmath}$ is a real $2 \times 2$ matrix of the form 
$\left[
\begin{array}{lr}
         a_{e+t+\jmath} & -b_{e+t+\jmath} \\
         b_{e+t+\jmath} & a_{e+t+\jmath}
        \end{array} \right] $ 
for $1 \le  \jmath \le f-t $. 
Notice that $\phi$ and $\psi$ have the same minimal polynomial over $\Q$, since $\phi$ is the diagonal matrix containing $J$ copies of $\psi$.
Let us consider now the following algebraic embeddings:
\begin{eqnarray*} \label{}
\lefteqn{\Psi_j: \Z[\phi]{\balpha_j}   \to  \mathbb{G}_j,} && \\
&& P_j(\phi){\balpha_j}  \mapsto P_j(D_j) \bbeta_j \,,
\end{eqnarray*}
where $P_j(x)$ is a polynomial over $\Z$ and $\bbeta_j :=  (1, \cdots, 1) \in  \mathbb{G}_j$.
Note that 
\[ \Psi_j(\phi x) = D_j \Psi_j(x) \ \ \ \mbox{for any $x \in \Z[\phi]{\balpha_j} $}\,.\]


Now we can define a map 
\begin{eqnarray} \label{Psi-map}
\Psi \ : \ \Z[\phi]{\balpha}_1 + \cdots +  \Z[\phi]{\balpha}_J  & \ \ \to \ \  &  {\mathbb{G}}_1 \times \cdots \times 
{\mathbb{G}}_J ,  \nonumber \\
P_1(\phi)\balpha_1 + \cdots + P_J(\phi)\balpha_J & \ \  \mapsto \ \ & (P_1(D_1) \bbeta_1, \dots, P_J(D_J) \bbeta_J).
\end{eqnarray}
Since $\balpha_1, \ldots, \balpha_J$ are linearly independent over $\Z[\phi]$, the map $\Psi$ is well defined.
Thus ${\Psi(\phi x) = D \Psi(x)}$ for 
\begin{eqnarray} \label{matrix-on-internalSpace}
D := \left[ 
\begin{array}{ccc} 
D_1 & \cdots & O \\
\vdots & \ddots & \vdots \\
O & \cdots & D_J \
\end{array}
\right],
\end{eqnarray}
where $D_1 = \cdots = D_J$. Let $\mathbb{G} := {\mathbb{G}}_1 \times \cdots \times 
{\mathbb{G}}_J $.

Let us construct a cut and project scheme :
\be \label{cut-and-project3}
\begin{array}{ccccc}
 \R^d & \stackrel{\pi_{1}}{\longleftarrow} & \R^d \times \mathbb{G} & \stackrel{\pi_{2}} {\longrightarrow} & \mathbb{G} \\
  && \cup &&\\
 L & \longleftarrow & {\widetilde{L}} & \longrightarrow & \Psi(L) \\
  & & & & \\
x & \longmapsfrom & (x, \Psi(x)) & \longmapsto & \Psi(x) \,,
\end{array}
\ee
where $\pi_1$ and $\pi_2$ are canonical projections, 
\[ L= \Z[\phi]{\balpha}_1 + \cdots +  \Z[\phi]{\balpha}_J \] and 
\[ \widetilde{L} = \{(x, \Psi(x)) : x \in L \}. \]
It is easy to see that $\pi_1|_{\widetilde{L}}$ is injective. 
We shall show that
$\pi_2(\widetilde{L})$ is dense in $\mathbb{G}$ and $\widetilde{L} $ is a lattice in $\R^d \times \mathbb{G}$.
We note that $\pi_2|_{\widetilde{L}}$ is injective, since $\Psi$ is injective.
Since $\phi$ commutes with the isomorphism $\sigma$ in Theorem \ref{controlPoints-Lattice}, 
we may identify $\mathcal{C}(\Tk)$ and its isomorphic image.
Thus from Theorem \ref{controlPoints-Lattice}, 
\[  \mathcal{C}(\Tk) = \bigcup_{i \le \kappa} \mathcal{C}_i \ \subset \ \Z[\phi]{\balpha}_1 + \cdots +  \Z[\phi]{\balpha}_J,    \]
where $\balpha_1, \dots, \balpha_J \in \mathcal{C}(\Tk)$.  
Note that for any $k \in \N$ and $1 \le j \le J$, $\phi^k \balpha_j \in \mathcal{C}(\Tk)$.
So we can note that 
\begin{eqnarray}  \label{C-generate-L}
L= \langle \bigcup_{i\le \kappa} \mathcal{C}_i \rangle.
\end{eqnarray}

\begin{lemma}
\label{lattice}
 $\widetilde{L}$ is a lattice in $\R^d \times \mathbb{G}$.
\end{lemma}

\proof

By Cayley-Hamilton theorem, there exists a monic polynomial $p(x)\in \Z[x]$ of degree $n$
such that $p(\phi)=id$. Thus every element of $\Z[\phi]$ is expressed as a polynomial of
$\phi$ of degree $n-1$ with integer coefficients
where the contant term is identified with a constant multiple of the identity matrix.
Therefore $L$ is a free $\Z$-module of rank $nJ$. 
Notice that $L$ and $\widetilde{L}$ are isomorphic $\Z$-modules so that $\widetilde{L}$ is also a free $\Z$-module of rank $nJ$ on $\R^d \times \mathbb{G}$. 
 Let us define 
 \[ \bgamma_{j} := (\balpha_j, \Psi(\balpha_j)) \in \R^d \times \mathbb{G} \ \ \ \mbox{for any $1 \le j \le J$}. \]
 Then, in fact, for any $1 \le k \le (e+f)J$,
 \be \label{def-beta}
 (\bgamma_{j})_{k} =
\left\{\begin{array}{ll}
                             1 \ \ \ & \mbox{if} \ \  (j-1)m + 1 \le k \le jm;\\
                             1 \ \ \ & \mbox{if} \ \  d+(j-1)(n-m) + 1 \le k \le d+j(n-m); \\ 
                             0 \ \ \ & \mbox{else} \,.
                           \end{array} \right.
\ee
Define also
 \[ \Omega = \left[ \begin{array}{cc}
 \phi & 0 \\
 0 & D
 \end{array} \right]
 \]
which is a linear map on $\R^d \times \mathbb{G}$. 
Note that $(\pi_1|_{\widetilde{L}})^{-1} (\Z[\phi]{\balpha}_j) = \Z[\Omega]\bgamma_j $ and  
 $\Z[\Omega]{\bgamma}_j$ is isomorphic to the image of $\Z^n$ by multiplication of the 
$n \times n$ matrix $A=(\lambda_i^{k-1})_{i,k \in \{1,\dots,n\}}$. 
Since $A$ is non-degenerate by the
Vandermonde determinant, 
\[ \{{\bgamma_1}, \dots, \Omega^{n-1}{\bgamma_1}, \dots, {\bgamma_J}, \dots, \Omega^{n-1}{\bgamma_J}\} \] 
forms a basis of $\R^d \times \mathbb{G}$ over $\R$. Thus $\widetilde{L}$ is a lattice in $\R^d \times \mathbb{G}$. 


\qed

\begin{lemma} \label{phi(L)-dense-in-K}
$\Psi(L) = \pi_2(\widetilde{L})$ and $\pi_2(\widetilde{L})$
is dense in $\mathbb{G}$.
\end{lemma}

\proof
For simplicity, we prove the totally real case, i.e., $\lam_i\in \R$ for all $i$. 
Since the diagonal blocks of $\phi$ are all same, it is enough to show that   
$\Psi_1(\Z[\phi] \balpha_1)$ is dense in $\mathbb{G}_1$. 
By (Theorem 24 (Siegel 1989)), $\Psi_1(\Z[\phi] \balpha_1)$ is dense in $\mathbb{G}_1$ if 
$$
\sum_{i=m+1}^{n} x_{i}
 \lam_i^{k-1}  \in \Z \quad (k=1, \dots, n)
$$
implies $x_{i}=0$ for $i=m+1,\dots, n$. The condition is equivalent to
$$
\xi A \in \Z^n
$$
with $\xi=(x_i)=(0,\dots,0,x_{m+1},\dots, x_n)\in \R^n $
in the terminology of Lemma \ref{lattice}. Multiplying by the inverse of $A$,
we see that the entries of $\xi$ must be Galois conjugates.
As $\xi$ has at least one zero entry, we obtain $\xi=0$ which shows
$x_i=0$ for $i=m+1,\dots, n$.
In fact, this discussion is using
the Pontryagin duality that the $\Theta : \Z^n \rightarrow \R^{n-m}$
has a dense image if and only if
its dual map $\widehat{\Theta}: \R^{n-m} \rightarrow \T^n$ is injective
(see also (Chapter II, Section 1 (Meyer 1972)), (Iizuka \& Akiyama \& Akazawa 2009)
and (Akiyama 1999)).
The case with complex conjugates is similar.\qed

\medskip

Now that we constructed the cut-and-project scheme (\ref{cut-and-project3}), we would like to introduce a special projected set $E_{\delta}$ which will appear in the proofs of the main results in Section \ref{MainResult}.
For $\delta > 0$, we define
\begin{eqnarray}  \label{projectedSetFromBall}
E_{\delta} &:=& \Lambda(B_{\delta}^{\mathbb{G}}(0))   \nonumber \\
&=& \pi_1({\pi_2}^{-1}(\Psi(L) \cap B_{\delta}^{\mathbb{G}}(0))) \nonumber \\
 &=& \{ \mathcal{P(\phi)  \balpha} \in  L : \Psi(\mathcal{P(\phi)  \balpha} ) \in B_{\delta}^{\mathbb{G}}(0) \}.
\end{eqnarray} 

In the following lemma, we find an adequate window for a set $\phi^n E_{\delta} $ and note that $E_{\delta}$ is a Meyer set.

\begin{lemma}
For any $\delta > 0$ and $n \in \N$, if $E_{\delta} = \Lambda(B_{\delta}^{\mathbb{G}}(0))$, then 
\[ \phi^n E_{\delta} 
 =  \{ \mathcal{R(\phi)  \balpha} \in  L : \Psi(\mathcal{R(\phi)  \balpha} ) \in D^n \cdot B_{\delta}^{\mathbb{G}}(0) \} \]
 and $E_{\delta}$ forms a Meyer set.
\end{lemma}

\proof    Note that
\begin{eqnarray}
\lefteqn{\Psi(\mathcal{P(\phi)  \balpha} )  \in B_{\delta}^{\mathbb{G}}(0) } \nonumber \\
& \Longleftrightarrow & D \cdot \Psi(\mathcal{P(\phi)  \balpha} )  \in  D \cdot B_{\delta}^{\mathbb{G}}(0)   \nonumber \\
& \Longleftrightarrow & \Psi( \phi (\mathcal{P(\phi)  \balpha}))  \in  D \cdot B_{\delta}^{\mathbb{G}}(0). 
 \end{eqnarray}
Notice that if $\phi$ is unimodular, then $\phi L = L = \phi^{-1} L$ and $\phi^{-1} \in \Z[\phi]$. 
Thus 
\begin{eqnarray} \label{phi-E_delta}
\ \ \ \ \  \phi E_{\delta} 
& = &  \phi \{ \mathcal{P(\phi)  \balpha} \in  L : \Psi(\mathcal{P(\phi)  \balpha} ) \in B_{\delta}^{\mathbb{G}}(0) \}    \nonumber \\
& = &  \{ \phi \mathcal{P(\phi)  \balpha} \in  L : \Psi(\phi \mathcal{P(\phi)  \balpha} ) \in D \cdot B_{\delta}^{\mathbb{G}}(0) \}   \label{phi-E_delta-1} \\
& = & \{ \mathcal{Q(\phi)  \balpha} \in  L : \Psi(\mathcal{Q(\phi)  \balpha}) \in D \cdot B_{\delta}^{\mathbb{G}}(0) \}. \label{phi-E_delta-2}
\end{eqnarray}
It is easy to see that the set in (\ref{phi-E_delta-1}) is contained in the set in (\ref{phi-E_delta-2}). The inclusion for the other direction is due to the fact that 
$\phi L = L$ and $\phi^{-1} \in \Z[\phi]$.
Hence for any $n \in \N$,
\begin{eqnarray} \label{E-delta-and-ball}
\ \ \ \ \  \phi^n E_{\delta} 
& = & \{ \mathcal{R(\phi)  \balpha} \in  L : \Psi(\mathcal{R(\phi)  \balpha} ) \in D^n \cdot B_{\delta}^{\mathbb{G}}(0) \}.  
\end{eqnarray}
Since (\ref{cut-and-project3}) is a cut-and-project scheme and $B_{\epsilon}^{\mathbb{G}}(0)$ is bounded, $E_{\delta}$ forms a Meyer set  for each $\delta > 0$ 
(see (Moody 1997)).
\qed

\section{Pure discrete spectrum, Meyer set, and Pisot family}  \label{KnownResults}

\begin{lemma} [Lemma 4.10 (Lee \& Solomyak 2008)] \label{pdd-imply-relDenseEV}
Let $\Tk$ be a tiling on $\R^d$. Suppose that $(X_{\Tk}, \R^d, \mu)$ has pure discrete dynamical spectrum. Then the eigenvalues for the dynamical system $(X_{\Tk}, \R^d, \mu)$ span $\R^d$. 
\end{lemma}

\begin{prop} [Proposition 6.6 (Lee \& Solomyak 2018)]  \label{relDense-Meyer}
Let $\Tk$ be a primitive substitution tiling on $\R^d$ with an expansion map
$\phi$. Suppose that all the eigenvalues of $\phi$ are algebraic integers. Assume that the set of eigenvalues of $(X_{\Tk}, \R^d, \mu)$ is relatively dense.
Then  $\mathcal{C}(\Tk)$ is a Meyer set.
\end{prop}

We note that ``repetitivity'' is not necessary for Proposition \ref{relDense-Meyer}.
Under the assumption that  $\Tk$ is a primitive substitution tiling on $\R^d$, the following implication holds;
\[
\mbox{Pure discrete spectrum} \ \Longrightarrow \ \mbox{Relative dense eigenvalues} \ \Longrightarrow \ \mbox{Meyer set} \ \Longrightarrow \ \mbox{FLC}  
\] 



\begin{defi}
{\em A set of
algebraic integers $\Theta = \{\theta_1, \cdots, \theta_r \}$ is a
{\em Pisot family} if for any $1 \le j \le r$, every Galois
conjugate $\gamma$ of $\theta_j$, with $|\gamma| \ge 1$, is
contained in $\Theta$. For $r=1$, with $\theta_1$ real and $|\theta_1|>1$, this reduces to $|\theta_1|$ being a real Pisot number, and for 
$r=2$, with $\theta_1$ non-real and $|\theta_1|>1$, to $\theta_1$ being a complex Pisot number. 
}
\end{defi}

Under the assumption of rigidity of $\Tk$, we can derive the following proposition from (Lemma 5.1 [Lee \& Solomyak 2012]) without additionally assuming repetitivity and FLC.

\begin{prop}  [Lemma 5.1 (Lee \& Solomyak 2012)] \label{pdd-pisotFamily}
Let $\Tk$ be a primitive substitution tiling on $\R^d$ with a diagonalizable expansion map $\phi$. Suppose that all the eigenvalues of $\phi$ are algebraic conjugates with the same multiplicity and $\Tk$ is rigid.
Then if the set of eigenvalues of $(X_{\Tk}, \R^d, \mu)$ is relatively dense, then the set of eigenvalues of $\phi$ forms a Pisot family.
\end{prop}

\section{Main result} \label{MainResult}

We consider a primitive substitution tiling on $\R^d$ with a diagonalizable expansion map $\phi$.
Suppose that all the eigenvalues of $\phi$ are algebraically conjugate with the same multiplicity $J$ and $\Tk$ is rigid. Additionally we assume that there exists at least one algebraic conjugate $\lambda$ of eigenvalues of $\phi$ for which $|\lambda| < 1$. Recall that 
\[ \Xi =  \bigcup_{i=1}^{\kappa} (\mathcal{C}_i - \mathcal{C}_i),  \]
where $\mathcal{C}_i$ is the set of control points of tiles of type $i$ and $1 \le i \le \kappa$.
By the choice of the control point set in (\ref{control-point-set-translation-set}), we note that $L = \langle \Xi \rangle$.

\begin{lemma}  \label{Xi-in-E_delta}
Assume that the set of eigenvalues of $\phi$ is a Pisot family. Then
$\Xi \subset E_{\delta}$ for some $\delta > 0$, where $E_{\delta}$ is given in (\ref{projectedSetFromBall}).
\end{lemma}

\proof
Since we are interested in $\Xi$ which is a collection of translation vectors, the choice of control point set $\mathcal{C}(\Tk)$ doesn't really matter. So we use the tile map (\ref{the-tile-map}) which sends a tile to the same type of tiles in $\Tk$. From (Lemma 4.5 (Lee \& Solomyak 2008)), for any $y \in \Xi$, 
\[ y = \sum_{n=0}^N \phi^n x_n, \ \ \ \mbox{ where $x_n \in U$ and $U$ is a finite subset in $L$}. \] 
Since $\phi$ is an expansive map and satisfies the Pisot family condition, the maps $\Psi_j$ and $\Psi$ are defined with all the algebraic conjugates of eigenvalues of $\phi$
whose absolute values are less than $1$. 
Thus $\Psi(\Xi) \subset B_{\delta}^{\mathbb{G}}(0)$ for some $\delta > 0$. 
From the definition of $E_{\delta}$ in (\ref{projectedSetFromBall}), $\Xi \subset E_{\delta}$.
\qed

\begin{lemma} \label{latticePoints-in-translationVectorSet}
Assume that $\Tk$ has pure discrete spectrum.
Then for any $y \in \langle \Xi \rangle$, there exists $\ell = \ell(y) \in \N$ such that $\phi^{\ell}y \in \Xi$.
\end{lemma}

\proof
Note from (\ref{C-generate-L}) that for any $k \in \Z_{\ge 0}$ and $\balpha_j \in \{ \balpha_1, \dots, \balpha_J \}$, $\phi^k \balpha_j$  is contained in $\Xi$.
Recall that $\langle \Xi \rangle \subset \Z[\phi]{\balpha}_1 + \cdots +  \Z[\phi]{\balpha}_J $, where ${\balpha}_1, \dots , {\balpha}_J \in \mathcal{C}(\Tk) $. So any element $y \in \langle \Xi \rangle$ is a linear combination of $\balpha_1, \phi \balpha_1, \dots, \phi^{n-1} \balpha_1, \dots, \balpha_J, \phi \balpha_J, \dots, \phi^{n-1} \balpha_J$ over $\Z$. 
Applying (\ref{translation-vectorSet-relation}) many times if necessary, we get that  for any $y \in \langle \Xi \rangle$, $\phi^{\ell} y \in \Xi$ for some $\ell = \ell(y) \in \N$. 
\qed

\begin{prop} \label{ppd-modelSet}
Let $\Tk$ be a primitive substitution tiling on $\R^d$ with an expansion map $\phi$. 
Under the assumption of the existence of CPS (\ref{cut-and-project3}),
if $\Tk$ has pure discrete spectrum, then there exists $K \in \N$ such that 
\begin{eqnarray} \label{E_delta-contained-Xi}
\phi^K E_{\delta} \subset \Xi. 
\end{eqnarray}
\end{prop}

\proof
We first prove that there exists a finite set $F$ such that for all $x \in E_{\delta}$, $x \in \Xi - v$ for some $v \in F$. This can be obtained directly from (Lemma 5.5.1 (Strungaru 2017)), but for reader's convenience we give the proof here. 
Note that $E_{\delta}$ is a Meyer set and $\Xi \subset E_{\delta}$ for some $\delta > 0$.
Since $\Xi$ is relatively dense, for any $x \in E_{\delta}$, there exists $r > 0$ such that $\Xi \cap B_r^{\R^d}(x) \neq \emptyset$.
From the Meyer property of $E_{\delta}$, the point set configurations  
\begin{eqnarray} \label{finite-configurations} 
 \{ \Xi \cap B_r^{\R^d}(x)  : x \in E_{\delta} \} 
 \end{eqnarray}
 are finite up to translation elements of $E_{\delta}$. We should note that if $E_{\delta}$ has FLC but not the Meyer property, the property (\ref{finite-configurations}) may not hold. 
Let
\begin{eqnarray*} 
F & = & \{ u - x : u \in \Xi \cap B_r^{\R^d}(x) \ \mbox{and} \ x \in E_{\delta}. \} 
\end{eqnarray*} 
Then 
\[ F = (\Xi - E_{\delta}) \cap B_r(0), \] 
$F \subset L$, and $F$ is a finite set.
Thus for any $x \in E_{\delta}$, 
\begin{eqnarray} \label{Meyer-implication} 
\mbox{$x \in \Xi - v$ \ for some $v \in F$}.
\end{eqnarray}

From Lemma \ref{latticePoints-in-translationVectorSet} and $L = \langle \Xi \rangle$, for any $y \in L$, there exists $\ell = \ell(y) \in N$ such that $\phi^{\ell}y \in \Xi$.
By the pure discrete spectrum of $\Tk$ and (\ref{translation-vectorSet-relation}), there exists $M \in \N$ such that 
\begin{eqnarray} \label{ppd-property1}
\phi^M \Xi - \phi^M \Xi \subset \Xi.
\end{eqnarray}
Applying the containment (\ref{ppd-property1}) finitely many times, we obtain that there exists $K_0 \in \N$ such that $\phi^{K_0} F \subset \Xi$.
Hence together with (\ref{Meyer-implication}), there exists $K \in \N$ such that 
\begin{eqnarray} \label{E_delta-contained-Xi}
\phi^K E_{\delta} \subset \Xi \,. 
\end{eqnarray}
\qed

In order to discuss about model sets and compute the boundary measures of their windows for substitution tilings, we need to introduce $\kappa-$set substitutions for substitution Delone sets which represent the substitution tilings.

\begin{defi}
{\em 
For a substitution Delone $\kappa$-set $\Lb = (\Lambda_i)_{i \le \kappa}$ satisfying (\ref{def-subst-mul}), define a matrix $\Phi=(\Phi_{ij})_{i,j=1}^{\kappa}$ 
whose entries are finite (possibly empty) families of linear affine transformations on 
$\R^d$ given by
$$
\Phi_{ij} = \{f: x\mapsto \phi x+a \ | \ a\in \Dk_{ij}\}\,.
$$
We define $\Phi_{ij}(\Xk) := \bigcup_{f\in \Phi_{ij}} f(\Xk)$ for  $\Xk\subset \R^d$. For a $\kappa$-set $(\Xk_i)_{i\le \kappa}$ let
\begin{eqnarray} \label{set-equations}
\Phi\bigl((\Xk_i)_{i\le \kappa}\bigr) = \Bigl(\bigcup_{j=1}^{\kappa} \Phi_{ij}(\Xk_j)\Bigr)_{i\le \kappa}\,.
\end{eqnarray}
Thus $\Phi(\Lb)= \Lb$ by definition. We say that $\Phi$ is a {\em $\kappa$-set substitution}.
Let 
\[ S(\Phi) = (\mbox{card} \Phi_{ij})_{ij} \] 
be a {\em substitution matrix} corresponding to $\Phi$. This is analogous to the substitution matrix for a tile-substitution.
}
\end{defi}

Recall that there exists a finite generating set $\mathbf{P}$ such that 
\begin{eqnarray}  \label{generatingSet}  
\mathcal{C} = \mbox{lim}_{r \to \infty}  \Phi^r(\mathbf{P})  
\end{eqnarray}
from (Lagarias \& Wang 2003, Lee \& Moody \& Solomyak 2003). 
If the finite generating set $\mathbf{P}$ consists of a single element, we say that $\mathcal{C}$ is {\em generated from one point}. 
Since $\Psi(L)$ is dense in $\mathbb{G}$ by Lemma \ref{phi(L)-dense-in-K}, we have a unique extension of $\Phi$ to a $\kappa$-set substitution on $\mathbb{G}$ in the obvious way; 
if $f \in \Phi_{ij}$ for which $f: L  \to L$, $f(x) = \phi x + a$, we define $f^*: \Psi(L) \to \Psi(L)$, $ f^*(u) = D u + a^* $, $D$ is given in (\ref{matrix-on-internalSpace}), and $a^* = \Psi(a)$. Since $\Psi(L)$ is dense in $\mathbb{G}$, we can extend the mapping $f^*$ to $\mathbb{G}$. If there is no confusion, we will use the same notation $f^*$ for the extended map.

Note that, by the Pisot family condition on $\phi$, there exists some $c < 1$ such that $|D x| \le c \cdot |x|$ for any $x \in \Psi(L)$. 
This formula defines a mapping on $\mathbb{G}$ and $f^*$ is a contraction on $\mathbb{G}$. Thus 
a $\kappa$-set substitution $\Phi$ determines a multi-component iterated function system $\Phi^*$ on $\mathbb{G}$.
Let $S(\Phi^*) = (\mbox{card}   (\Phi_{ij}^*))_{ij}$ be a {\em substitution matrix} corresponding to $\Phi^*$.
 Defining the compact subsets 
\[  V_i = \overline{\Psi(\mathcal{C}_i)}  \ \ \ \mbox{for each $1 \le i \le \kappa$ } \]
and using (\ref{set-equations}) and the continuity of the mappings, we have 
\[  V_i =  \bigcup_{j=1}^{\kappa} \bigcup_{f^* \in (\Phi^*)_{ij}} f^*(V_j), \ \ \ i = 1, \dots, \kappa . \]
This shows that $V_1, \dots, V_{\kappa}$ are the unique attractor of $\Phi^*$.

\medskip

\begin{remark} \label{algCoincidence-at-origin}
{\em From (Prop. 4.4 (Lee 2007)), if $\Tk$ has pure discrete spectrum, then there exists $\mathcal{R} \in X_{\Tk}$ such that 
the control point set $\mathcal{C}_{\mathcal{R}}:= \mathcal{C}(\mathcal{R})$ of the tiling $\mathcal{R}$ satisfies
\[ \mathcal{C}_{\mathcal{R}} = \lim_{n \to \infty} (\Phi^N)^n (y)   \ \ \ \mbox{and} \ \ \
 y + \phi^N \Xi(\mathcal{R}) \subset (\mathcal{C}_{\mathcal{R}})_j \]
for some $y \in (\mathcal{C}_{\mathcal{R}})_j$, $j \le \kappa$ and $N \in \Z_+$. 
Note that $\omega^N (\mathcal{R}) = \mathcal{R}$. Let $\upsilon = \omega^N$.
We can consider a $r$-th level supertiling $\upsilon^{r} (\mathcal{R})$ of $\mathcal{R}$. Note that there exists an $r$-th level supertile $\upsilon^{r}(S)$ in $\upsilon^{r} (\mathcal{R})$ containing the origin in the support which contain the tile $y + T_j \in \mathcal{R}$. Redefining the tile map for the control points of this supertiling so that the control point of $r$-th level supertile $\upsilon^{r}(S)$ is at the origin, we can build a substitution tiling $\mathcal{R}' \in X_{\Tk}$ for which algebraic coincidence occurs at the origin.  
So rewriting the substitution if necessary, we can assume that $y=0$.
With this assumption, we get the following proposition.
}
\end{remark}

\begin{prop} \label{modelSetWithOpenWindow}
Let $\Tk$ be a primitive substitution tiling on $\R^d$ with a diagonalizable expansion map $\phi$ which is unimodular. Suppose that all the eigenvalues of $\phi$ are algebraic conjugates with the same multiplicity and $\Tk$ is rigid.
Suppose that 
\begin{eqnarray} \label{projectedPointSet-from-OpenSet} 
\mathcal{C} = \lim_{n \to \infty} (\Phi^N)^n (\{0\})   \ \ \ \mbox{and} \ \ \
  \phi^N \Xi \subset \mathcal{C}_j 
\end{eqnarray}
for some $0 \in \mathcal{C}_j$, $j \le \kappa$ and $N \in \Z_+$. 
Assume that CPS (\ref{cut-and-project3}) exists. Then each point set 
\begin{eqnarray} \label{ProjectedPointSet}
\mathcal{C}_i = \Lambda(U_i), \ \ \ i \le \kappa
\end{eqnarray}
 is an Euclidean model set in CPS (\ref{cut-and-project3}) with a  window $U_i$ in $ \mathbb{G}$ which is open and precompact.
\end{prop}

\proof
For each $i \le \kappa$ and $z \in \mathcal{C}_i$, there exists $n \in \Z_+$ such that 
\[ f(0) = z \ \ \ \mbox{and} \ \ f \in (\Phi^n)_{ij}. \] 
From $\phi^N \Xi \subset  \mathcal{C}_j $, 
\[ z + \phi^{n+N} \Xi \subset \mathcal{C}_i .\]
By Theorem \ref{pdd-algCoin} and Proposition \ref{ppd-modelSet}, there exists $K \in \N$ such that 
$\phi^K E_{\delta} \subset \Xi$. 
Thus 
\[ \mathcal{C}_i  = \bigcup_{z \in  \mathcal{C}_i} (z+ \phi^{N_z} E_{\delta_z}), \]
where $N_z$ depends on $z$.
From the equality of (\ref{E-delta-and-ball}), we let 
\[ U_i := \bigcup_{z \in  \mathcal{C}_i} (z^* + D^{N_z} B_{\delta_z}(0))  \ \ \ 
\mbox{for any $i \le \kappa$}.  \] 
Then 
\[  \mathcal{C}_i  =  \Lambda(U_i) \ \ \ \mbox{where $U_i$ is an open set in $\mathbb{G}$}, \]
for any $i \le \kappa$.

From Lemma \ref{Xi-in-E_delta}, $\Xi \subset E_{\delta}$ for some $\delta > 0$.
Thus $\overline{\Psi(\Xi)} \subset \overline{B_{\delta}^{\mathbb{G}}(0)}$. 
Since $ \overline{B_{\delta}^{\mathbb{G}}(0)}$ is compact, $\overline{\Psi(\Xi)}$ is compact. 
Thus $\overline{\Psi(\mathcal{C}_i)}$ is compact. 
\qed

\medskip

We can assume that the open window $U_i$ in (\ref{ProjectedPointSet}) is the maximal element satisfying (\ref{ProjectedPointSet}) for the purpose of proving the following proposition.

\begin{prop} \label{ppd-regularModelSet}
Let $\Tk$ be a repetitive primitive substitution tiling on $\R^d$ with a diagonalizable expansion map $\phi$ which is unimodular. Suppose that all the eigenvalues of $\phi$ are algebraic conjugates with the same multiplicity and $\Tk$ is rigid.
Under the assumption of the existence of CPS (\ref{cut-and-project3}), if 
\begin{eqnarray} \label{projectedPointSet-from-OpenSet} 
\mathcal{C} = \lim_{n \to \infty} (\Phi^N)^n (\{0\})   \ \ \ \mbox{and} \ \ \
 \phi^N \Xi \subset \mathcal{C}_j 
\end{eqnarray}
where $0 \in \mathcal{C}_j$, $j \le \kappa$ and $N \in \Z_+$,
then each Euclidean model set $\mathcal{C}_j$, $1 \le j \le \kappa$ has a window with boundary measure zero in the Euclidean internal space $\mathbb{G}$ of CPS (\ref{cut-and-project3}).
\end{prop}

\proof
Let us define $W_i = \overline{U_i}$,
where $U_i$ is the maximal open set in $\mathbb{G}$ satisfying (\ref{ProjectedPointSet}).
From the assumption of (\ref{projectedPointSet-from-OpenSet}), we first note that $\phi$ fulfills the Pisot family condition from Theorem \ref{pdd-algCoin} and Proposition \ref{pdd-pisotFamily}. 
For every measurable set $E \subset \mathbb{G}$ and for any $f^* \in (\Phi^*)_{ij}$ with $f^*(u) = Du+a^*$,
\[ \mu(f^*(E)) = \mu(D(E) + a^*) = |\mbox{det} D| \mu(E), \]
where $\mu$ is a Haar measure in $\mathbb{G}$ and $D$ is the contraction as given in (\ref{matrix-on-internalSpace}). 
Note that $|\mbox{det} D| < 1$.
In particular, 
\[ \mu(f^*(W_j)) = |\mbox{det} D| \mu(W_j), \ \ \ 1 \le j \le \kappa . \]
We have attractors $W_j$'s satisfying
\[ W_i = \bigcup_{j=1}^{\kappa} \bigcup_{f^* \in {(\Phi^*)}_{ij}} f^*(W_j). \]
Let us denote $w_j = \mu(W_j)$ for $1 \le j \le \kappa$ and $\mathbf{w} = [w_1, \dots, w_{\kappa}]^T$.
Then for any $r \in \N$,
\[  w_i \le \sum_{j=1}^{\kappa} |\mbox{det} D|^r \mbox{card}({(\Phi^*)}^r)_{ij} w_j. \] 
Note here that for any $1 \le j \le \kappa$, $w_j > 0$ follows from the fact that $W_j$ has non-empty interioir.
Thus 
\[ \mathbf{w} \le |\mbox{det} D|^r S({(\Phi^*)}^r) \mathbf{w} \le |\mbox{det} D|^r  (S({\Phi^*}))^r \mathbf{w} \ \ \ \mbox{for any $r \in \N $}.
\]
Note from (Lagarias \& Wang 2003) that Perron eigenvalue of $(S({\Phi^*}))^r$ is $|\mbox{det} \phi|^r$.  
From the unimodular condition of $\phi$, 
\[ \mbox{det} D \cdot  \mbox{det} \phi   = \pm 1 .  \]
Since $(S({\Phi^*}))^r$ is primitive, from (Lemma 1 (Lee \& Moody 2001))
\[ \mathbf{w} =  |\mbox{det} D|^r S({(\Phi^*)}^r) \mathbf{w} =  |\mbox{det} D|^r (S({\Phi^*}))^r \mathbf{w} \ \ \  \mbox{for any 
$r \in \N$}. 
\]
By the positivity of $\mathbf{w}$ and $S({(\Phi^*)}^r) \le (S({\Phi^*}))^r$, $S({(\Phi^*)}^r) = (S({\Phi^*}))^r$.

Recall that for any $r \in \N$, 
\begin{eqnarray}  \label{r-th-iterations-windows}
 W_i = \bigcup_{j=1}^{\kappa} ({(\Phi^*)}^r)_{ij} W_j . 
\end{eqnarray}
From (\ref{M-th-iteration-digitSet}), for any $r \in \N$,
\[  (\mathcal{D}^r)_{ij} = \bigcup_{k_1, k_2, \dots, k_{(r-1)} \le \kappa} (\mathcal{D}_{ik_1} + \phi \mathcal{D}_{k_1 k_2} +
\cdots + \phi^{r-1} \mathcal{D}_{k_{r-1}j}) 
\]
and
\[ (\Phi^r)_{ij}(\{x_j\}) = \phi^r x_j + (\mathcal{D}^r)_{ij}  \ \ \  \mbox{for any $x_j \in \mathcal{C}_j$}. \]
Note that $W_i = \overline{U_i} = \overline{\Psi(\mathcal{C}_i)}$ and $U_i$ is a non-empty open set. 
As $r \to \infty$, $\bigcup_{j=1}^{\kappa} (\mathcal{D}^r)_{ij}$ is dense in $W_i$.
Since $\mathbb{G}$ is an Euclidean space,  we can find a non-empty open set $V \subset \mathbb{G}$ such that 
$V \subset \overline{V} \subset U_i$.
So there exist $M \in \N$  and $a \in (\mathcal{D}^M)_{ij}$ such that $a^* + D^M \Psi(\Xi)  \subset V \subset \overline{V} \subset U_i$.
Since $ W_j \subset \overline{\Psi(\Xi)}$, 
\[ a^* + D^M(W_j) \subset a^* + D^M \overline{\Psi(\Xi)}. \]
Thus there exists $g^* \in ((\Phi^*)^M)_{ij}$ such that  
\begin{eqnarray} \label{someWindow-in-otherWindow} 
g^*(W_{j}) \subset \ U_i. 
\end{eqnarray}
Hence
\begin{eqnarray} 
\partial{U_i} = W_i \backslash U_i & = & \bigcup_{j=1}^{\kappa} (({\Phi^*})^M)_{ij} (W_j) \ \backslash \ U_i  \nonumber \\
& \subset & \bigcup_{j=1}^{\kappa} \left( (({\Phi^*})^M)_{ij} (W_j) \right) \ \backslash \  \left(  (({\Phi^*})^M)_{ij} (U_j)  \right) 
\label{second-inclusion} \\
& \subset & \bigcup_{j=1}^{\kappa} (({\Phi^*})^M)_{ij} (\partial{U_j}). \label{boundary-containment}
\end{eqnarray}
The inclusion (\ref{second-inclusion}) is followed by the maximal choice of an open set $U_i$.
Let 
\[ v_j = \mu(\partial U_j), \ 1 \le j \le \kappa, \ \  \mbox{and} \ \ \mathbf{v} = [v_1, \dots, v_{\kappa}]^T,  . \]
Then \[ \mathbf{v} \le  |\mbox{det} D|^M S({\Phi^*})^M \mathbf{v} . \] 
From (\ref{someWindow-in-otherWindow}), we observe that not all functions in $S({(\Phi^*)}^M)$ are used for the inclusion (\ref{boundary-containment}).
Thus there exists a matrix $S'$ for which 
\[ 0 \le \mathbf{v} \le  |\mbox{det} D|^M S' \mathbf{v} \le   |\mbox{det} D|^M S({(\Phi^*)}^M) \mathbf{v}  \le  |\mbox{det} D|^M (S({\Phi^*}))^M \mathbf{v} \,, \]
where $S' \le (S({\Phi^*}))^M $ and $S' \neq (S({\Phi^*}))^M $.
If $\mathbf{v} > 0$, again from (Lemma 1 (Lee \& Moody 2001)), $S' = (S({\Phi^*}))^M $. This is a contradiction to (\ref{someWindow-in-otherWindow}).
Therefore $v_j= 0 $ for any ${1 \le j \le \kappa}$. 
\qed

\medskip

The regularity property of model sets can be shared for all the elements in $X_{\Tk}$. One can find the earliest result of this property in (Schlottmann 2000) and the further development in (Baake \& Lenz \& Moody 2007, Lee \& Moody 2006). We state the property (Prop. 4.4 (Lee \& Moody 2006)) here.

\begin{prop} (Schlottmann 2000, Baake \& Lenz \& Moody 2007, Lee \& Moody 2006) \label{window-containment}
Let $\mathcal{C}$ be a Delone $\kappa$-set in $\R^d$ for which $\Lambda({V_i}^{\circ}) \subset \mathcal{C}_i \subset \Lambda(\overline{V_i})$ where  
$\overline{V_i}$ is compact and ${V_i}^{\circ} \neq \emptyset$ for $i \le \kappa$ with respect to to some CPS. Then for any $\Gb \in X_{\mathcal{C}}$, there exists $(-s, -h) \in \R^d \times \mathbb{G}$ so that
\[ -s + \Lambda(h + {V_i}^{\circ}) \subset \Gamma_i \subset -s + \Lambda(h + \overline{V_i}) \ \ \ \mbox{for each $i \le \kappa$.} \]
\end{prop}
\qed

\medskip

From the assumption of pure discrete spectrum and Remark \ref{algCoincidence-at-origin}, we can observe that the condition  (\ref{projectedPointSet-from-OpenSet}) is fulfilled in the following theorem.

\begin{theorem}  \label{ppd-regModelSet}
Let $\Tk$ be a repetitive primitive substitution tiling on $\R^d$ with a diagonalizable expansion map $\phi$ which is unimodular.
Suppose that all the eigenvalues of $\phi$ are algebraically conjugate with the same multiplicity.
If $\Tk$ has pure discrete spectrum, then each control point set $\mathcal{C}_j$, $1 \le j \le \kappa$ is a regular Euclidean model set in CPS (\ref{cut-and-project3}).  
\end{theorem}

\proof
Under the assumption of pure discrete spectrum, we know that $\Tk$ has FLC from (Lee \& Solomyak 2018) and $\phi$ fulfills the Pisot family condition from (Lee \& Solomyak 2012). 
From Theorem \ref{controlPoints-Lattice}, we know that $\Tk$ is rigid. 
Since $\phi$ is unimodular, there exists at least one algebraic conjugate $\lambda$ of eigenvalues of $\phi$ for which $|\lambda| < 1$.
Thus we can construct the CPS (\ref{cut-and-project3}) with an Euclidean internal space. Since $\Tk$ has pure discrete spectrum and is repetitive, we can find a substitution tiling $\mathcal{S}$ in $X_{\Tk}$ such that 
\begin{eqnarray} 
\mathcal{C}_{\mathcal{S}} = \lim_{n \to \infty} (\Phi^N)^n (\{0\})   \ \ \ \mbox{and} \ \ \
 \phi^N \Xi \subset {(\mathcal{C}_{\mathcal{S}})}_j 
\end{eqnarray}
where $0 \in ({\mathcal{C}_{\mathcal{S}})}_j$, $j \le \kappa$ and $N \in \Z_+$.
The claim follows from Proposition \ref{ppd-modelSet}, \ref{ppd-regularModelSet} and \ref{window-containment}.

\begin{cor}
Let $\Tk$ be a repetitive primitive substitution tiling on $\R^d$ with a diagonalizable expansion map $\phi$ which is unimodular.
Suppose that all the eigenvalues of $\phi$ are algebraically conjugate with the same multiplicity.
Then $\Tk$ has pure discrete spectrum if and only if each control point set $\mathcal{C}_j$, $1 \le j \le \kappa$ is a regular Euclidean model set in CPS (\ref{cut-and-project3}).  
\end{cor}

\proof
It is known that any regular model sets have pure discrete spectrum in quite a general setting (Schlottmann 2000). Together with Theorem \ref{ppd-regModelSet}, we obtain the equivalence between pure discrete spectrum and regular model set in substitution tilings.
\bigskip

The next example shows that the unimodularlity of $\phi$ is necessary.

\begin{exam} \label{non-unimodular}
{\em Let us consider an example of non-unimodular substitution tiling which is studied in (Baake \& Moody \& Schlottmann 1998). 
This example is proven to be a regular model set in the setting of a cut-and-project scheme constructed in (Baake \& Moody \& Schlottmann 1998) with the help of $2$-adic embedding.
In our setting of CPS (\ref{cut-and-project3}), we show that this example can not provide a model set,
since we are only interested in the Euclidean window in this paper. 

The substitution matrix of the primitive two letter substitution 
\[  a \to aab \ \ \  \ \    b \to abab   \]
has the Perron-Frobenius eigenvalue $\beta:=2+\sqrt{2}$ which is a Pisot number but non-unimodular.  
We can extend the letter $a$ to the right hand side by the substitution and the letter $b$ to the left hand side. So we can get a bi-infinite sequence fixed under the substitution.
A geometric substitution tiling arising from this substitution can be obtained by replacing symbols $a$ and $b$ in this sequence by the intervals of length $\ell(a) =1$ and $\ell(b) = \sqrt{2}$. Then we have the following tile-substitution 
$\omega$
\[   \omega(T_a) = \{ T_a, 1+ T_a, 2+ T_b \} , \]
\[   \omega(T_b) = \{T_a, 1+ T_b, 1+\sqrt{2}+T_a, 2+\sqrt{2} + T_b \} , \]
where $T_a = ([0,1], a)$ and $T_b=([0, \sqrt{2}], b)$.
Considering return words $\{a,ab\}$ for $a$, and $\{ba,baa\}$ for $b$, 
we can check $L=\langle \Xi\rangle$. 
We choose left end points $\Lam_a, \Lam_b$ of corresponding intervals as the set of control points.
Then they satisfy
$$
\Lam_a=\left(\beta \Lam_a +\{0,1\}\right) \bigcup \left(\beta \Lam_b +\{0,1+\sqrt{2}\}\right)
$$
$$
\Lam_b=\left(\beta \Lam_a +\{2\}\right) \bigcup \left(\beta \Lam_b+\{1,2+\sqrt{2}\}\right),
$$
by Lagarias-Wang duality (Lagarias \& Wang 2003). Applying Galois conjugate $\kappa$ which sends $\sqrt{2}\to -\sqrt{2}$, 
we obtain a generalized iterated function system
$$
X_a=\left(\gamma X_a +\{0,1\}\right) \bigcup \left(\gamma X_b +\{0,1-\sqrt{2}\}\right)
$$
$$
X_b=\left(\gamma X_a +\{2\}\right) \bigcup \left(\gamma X_b+\{1,2-\sqrt{2}\}\right),
$$
with $\gamma=2-\sqrt{2}$, $X_a=\overline{ \kappa(\Lam_a)}$ and $X_b=\overline{\kappa(\Lam_b)}$.
We can easily confirm that
$$
X_a=[0,1+\sqrt{2}], X_b=[\sqrt{2},2+\sqrt{2}]. 
$$
are the unique attractors of this iterated function system. 
Since $X_a\cap X_b=[\sqrt{2},1+\sqrt{2}]$
contains an inner point, it is unable to distinguish them by any window in this setting.
}
\end{exam}

\section{Further study}

\noindent
We have mainly considered unimodular substitution tilings in this paper. 
Example \ref{non-unimodular} shows a case of non-unimodular substitution tiling which is studied in (Baake \& Moody \& Schlottmann 1998).  
It cannot be an Euclidean model set in the cut-and-project scheme (\ref{cut-and-project3}) that we present in this paper,  
but it is proven to be a regular model set in the setting of a cut-and-project scheme constructed in (Baake \& Moody \& Schlottmann 1998),
which suggests non-unimodular tilings require non archimedian embeddings to construct internal spaces.
It is an intriguing open question to construct a concrete cut-and-project scheme in this case.

\medskip

\appendix

\subsection*{Acknowledgments}

We would like to thank to F. G\"ahler, U. Grimm, M. Baake, and N. Strungaru for the valuable and important comments and discussions at MATRIX in Melbourne. The third author would like to also thank to A. Quas for his interest and questions for this work at ESI in Austria. We are grateful to the MATRIX and the ESI for their hospitality. We are indebted to the two referees for their important comments which improve the paper much more.
The work of D.-i. Lee was supported by a research grant from Seoul Women's University (2020-0205).
S. Akiyama was partially supported by JSPS grants (17K05159, 17H02849, BBD30028).
The research by J.-Y. Lee was supported by NRF grant No. 2019R1I1A3A01060365.
She is grateful to the KIAS, where part of this work was done.

\vspace{6mm}

\end{document}